\documentclass{amsart}
\usepackage{amsmath,amsthm,amssymb}
\usepackage[dvips]{color}
\usepackage{amsfonts}
\usepackage{rotating}
\usepackage{euscript}
\usepackage{pst-node}
\usepackage{epsfig}

\def\be{\begin{equation}}
\def\ee{\end{equation}}

\def\tilde{\widetilde}

\def\C{{\mathbb C}} 
 
\def\P{{\mathbb P}}
\def\Z{{\mathbb Z}}
\def\R{{\mathbb R}} 
\def\Q{{\mathbb Q}}

\def\f{\EuScript}
\def\e{\eqref}
\def\phi{{\varphi}}
\def\v{{\varepsilon}} 
\def\tt{\widetilde}
\def\deg{{\rm deg\,}}

\def\exp{{\rm exp\,}}

\def\mod{{\rm mod\ }}

\def\bp{\begin{proposition}}
\def\ep{\end{proposition}}

\def\bt{\begin{theorem}}
\def\et{\end{theorem}}
\def\br{\begin{remark}}
\def\er{\end{remark}}
\def\be{\begin{equation}}
\def\bee{\begin{equation*}}
\def\l{\label}
\def\e{\eqref}

\def\ee{\end{equation}}
\def\eee{\end{equation*}}
\def\bl{\begin{lemma}}
\def\el{\end{lemma}}
\def\bc{\begin{corollary}}
\def\ec{\end{corollary}}
\def\pr{\noindent{\it Proof. }}

\def\bd{\begin{definition}}
\def\ed{\end{definition}}
\input epsf.sty

\newtheorem{theorem}{Theorem}[section]
\newtheorem{lemma}{Lemma}[section]
\newtheorem{definition}{Definition}[section]
\newtheorem{corollary}{Corollary}[section]
\newtheorem{proposition}{Proposition}[section]
\newtheorem{remark}{Remark}[section]

\begin{document}
\title{Moments on Riemann surfaces and  hyperelliptic Abelian integrals}
\author{L. Gavrilov}
\address{Institut de Math\'{e}matiques de Toulouse, UMR 5219\\
  Universit\'{e}  de Toulouse,  31062 Toulouse,  France}
%\curraddr{}
\email{lubomir.gavrilov@math.univ-toulouse.fr}
\thanks{Part of the paper was written while the first author was visiting the Ben-Gurion University. He is obliged for the hospitality}
\author{F. Pakovich}
\address{Department of Mathematics, Ben-Gurion University
of the Negev\\ P.O.B. 653, Beer-Sheva, Israel; e-mail:
%{\texttt pakovich@cs.bgu.ac.il}
}
%\curraddr{}
\email{pakovich@math.bgu.ac.il}
%\thanks{thank you}

%    author two information
\date{\today}
%\abstract{ab}
\input epsf.sty
\begin{abstract}
In the present paper we solve the following different but interrelated problems: (a) the moment problem on Riemann surfaces, (b) the vanishing problem of polynomial Abelian integrals of dimension zero on the projective plane, (c) the  vanishing problem of polynomial hyperelliptic Abelian integrals.
\end{abstract}
\maketitle
\tableofcontents
\section{Introduction}
In the present paper we solve the following different but interrelated problems: 
\begin{description}
\item[(a)]
Let $f$ and $\omega$ be a function and a one-form meromorphic on a compact Riemann
surface $R$ and $\gamma \subset R$
be a closed curve. In Section 2 we give necessary and sufficient
conditions for all the ``moments''
\be \l{zer} m_s=\int_{\gamma} f^s \omega , \ \ \ s \geq 0,\ee
to vanish identically.
These condition are expressed in terms of the vanishing of a finite collection of algebraic functions, which can be interpreted as Abelian integrals of dimension zero on $R$.

\item[(b)]
Motivated by problem {\bf(a)}, we describe in Section \ref{abint} necessary and sufficient conditions for the identical vanishing of  polynomial  Abelian integrals of dimension zero on the projective plane.
\item[(c)] Finally, in Section 4 we apply the results obtained to the problem of the identical vanishing of hyperelliptic Abelian integrals of the form
\begin{equation}
\label{abel1}
I(t)= \int_{\gamma(t)} P(x,y) dx + Q(x,y) dy, \quad P,Q \in \C[x,y]
\end{equation}
where  $\gamma(t) \in H_1(\Gamma_t,\mathbb{Z})$ is a continuous family of 1-cycles and
$$\Gamma_t= \{(x,y)\in \C^2: y^2-f(x)=t\}, \quad f\in \C[x]$$ is a family of hyperelliptic curves. 
\end{description}
The moment vanishing problem {\bf (a)} has been  studied by several authors in the last years,
in relation with the center problem for the Abel equation
$$
\frac{dy}{dz}= p(z)y^2+ q(z) y^3, \quad p,q \in \C[z],
$$
see \cite{bry10} for an extensive list of references.
The solution that we present here in terms of zero-dimensional Abelian integrals is inspired by the approach of \cite{bull}, \cite{mom}.

Zero-dimensional Abelian integrals were introduced recently in \cite{gamo07}, in an attempt to verify certain conjectures concerning the 16th Hilbert problem in dimension zero. In particular, the problem of identical vanishing of such zero-dimensional integrals on the Riemann sphere for simple cycles has been studied and solved in \cite{gamo07, krma10}.
Notice however that in this case the problem reduces to the finding of conditions implying that for a pair of polynomials 
$f,g$ the equality 
$g(f^{-1}_i(z))\equiv g(f^{-1}_j(z))$ holds for two different branches of the function inverse to $f$, and in such a form the problem was studied and solved earlier (see e.g. \cite{roy}, \cite{let}).
In the general case the problem {\bf (b)} is to find conditions implying the equality 
$$\sum_{i=1}^na_ig(f^{-1}_i(z))\equiv 0,
 \ \ \ n=\deg f,$$ for arbitrary given rational $a_i\in \Q$. Essentially, in an implicit form a solution of this problem 
was already done in \cite{pm} as an ingredient of the proof of the so called polynomial moment problem. However, 
having in mind possible applications, for the reader convenience we present here a detailed and full exposition which 
is self-contained up to 
a single purely algebraic result of \cite{pm}.  Notice that 
a partial solution of problem {\bf (b)}, also basing on results of \cite{pm}, has been recently done in \cite{stm11}.

The last problem {\bf (c)} we solve concerns the identical vanishing of complete hyperelliptic Abelian integrals of the form $(\ref{abel1})$. Although this problem is of independent interest, we are once again motivated by applications to the 16th Hilbert problem. Namely, it is well known that if a limit cycle of the perturbed plane foliation
\begin{equation}
\label{fol}
d(y^2-f(x)) + \varepsilon (P(x,y) dx + Q(x,y) dy)= 0, \quad \varepsilon \sim 0
\end{equation}
bifurcates from the period orbit $\gamma(t_0)\subset \Gamma_{t_0}$ of the non-perturbed foliation, then the Abelian integral $I(t)$ defined by  (\ref{abel1}) vanishes at $t_0$. This is a corollary of the representation
\begin{equation}
\label{ret}
P_\varepsilon(t)=t+ \varepsilon I(t) + o(\varepsilon)
\end{equation}
of the first return map $P_\varepsilon$ associated to the family of periodic orbits $\gamma(t)$. The situation in which $I(t)\equiv 0$ is  exceptional, and this phenomenon is related to the singularities of the algebraic set of plane integrable foliation. Notice that even for plane polynomial vector fields of degree three,  the structure of the set of integrable foliations is unknown \cite{mova04} (for quadratic vector fields it is due to Dulac (1908) and Kapteyn (1912)). On the other hand, the identical vanishing of $I(t)$ only shows that the foliation (\ref{fol}) is integrable ``at a first order", 
and the study of the higher order terms  in the expansion (\ref{ret}) is needed in order to solve the associated center problem on the plane \cite{gavr05,gmn09}. 

The key idea to solve problem {\bf (c)} is to interpret the derivatives of $I(t)$ as moments \eqref{zer} 
for a certain choice of $R,f$ and $\omega$. Then the identical vanishing of $I(t)$ turns out equivalent, according to {\bf (a)}, to the identical vanishing of a collection of Abelian integrals of dimension zero. Furthermore, these Abelian integrals essentially reduce to the ones studied in {\bf (b)}. 
%Thus,
%combined all together, our results provide a complete solution of the vanishing problem of complete hyperelliptic Abelian integrals $(\ref{abel1})$. 

\section{Moments on Riemann surfaces and zero-dimensional Abelian integrals}
\label{section2}
Let $f, \omega$ be respectively a function  and a one-form meromorphic on a compact Riemann surface $R$, and let $\gamma\subset R$ be a closed rectifiable curve which avoids the poles of $f$. 
Then the moments (\ref{zer})
are well defined. In this section we shall suppose, for simplicity, that
the set of poles of $\omega$ is contained in the set of poles of $f$.
We shall give necessary and sufficient 
conditions for the generating function 
\be 
\label{generating0}
J(t)= J(\omega,f,\gamma,t)= -\sum_{i=0}^{\infty}\frac{m_i}{t^{i+1}}= \int_{\gamma} \frac{\omega}{f-t},\ \  t\sim \infty
\ee
of the moments $m_i$
to vanish identically or, more generally, to be rational. For this purpose we follow closely \cite{mom}, where the genus zero case, $R=\mathbb{CP}^1$, was studied in details.

Consider the induced map 
\begin{equation}\label{fmap}
\begin{array}{rcl}
f: R & \rightarrow & \mathbb{CP}^1 \\
     x & \mapsto  & [f(x):1]
\end{array}
\end{equation}
and let $\{c_1,c_2,\dots c_k\}$ be the set of all \emph{finite} critical values of $f$. For a regular generic value $c_0\in \mathbb{C}$, consider the ``star" $S\subset \mathbb{C}$ consisting of the segments $[c_0,c_i]$, $i=1, 2\dots k$.
It is not hard to show, using the assumption that $S$ 
contains all finite critical values of $f$, that one can continuously deform the path $\gamma$, without 
changing the corresponding function $J(t)$,   
in such a way that the image $f(\gamma)$ will be 
contained in $S$ (the explicit construction is given below).
Therefore (\ref{zer}) becomes
\be
\label{ms2}
m_s=\int_{\gamma} f^s \omega = \int_{\gamma} f^s \frac{\omega}{df} df = \sum_{i=1}^k \int_{c_0}^{c_i} \varphi_i(z)z^s dz
\ee
where $\varphi_i$ is an appropriate sum of branches of the algebraic function $$\frac{\omega}{df} \circ f^{-1}$$ 
in some simply-connected domain $U$ containing 
$S\setminus 
\{c_1,c_2,\dots c_k\}.$

Clearly, 
\be
\label{generating}
J(t) = \sum_{i=1}^k J_i(t), \mbox{  where  }\;\; J_i(t)= \int_{c_0}^{c_i}\frac{\varphi_i(z)}{z-t} dz.
\ee
Further, the functions $J_i(t)$ and therefore $J(t)$ allow for an analytic continuation on $\mathbb{C}\P^1\setminus 
\{c_1,c_2,\dots c_k\}.$
On the other hand, by a well-known property of Cauchy type integrals, 
the limits of the function $J(t)$ when $t$ approaches 
to a point $t\in [c_0,c_i]$ from the ``left"
and ``right'' sides of $[c_0,c_i]$ are related by the equality
$$J_i^+(t)-J_i^-(t)= 2\pi\sqrt{-1}\, \varphi_i(t).$$
Therefore, if the generating function $J(t)$ vanishes identically (or just allows for a single-valued analytical continuation), then the algebraic functions $\varphi_i$, $i=1,2 \dots k$ defined by (\ref{ms2}) vanish identically. Of course, $\varphi_i=0$ also implies that $J=0$.  

The study of conditions implying the vanishing of the algebraic functions $\varphi_i$ is \emph{a priori} a simpler problem than the initial one.
Furthermore, the functions $\varphi_i$ allow for a remarkable interpretation as $0$-dimensional Abelian integrals which we describe now.

Consider the singular fibration (\ref{fmap}) with fibers 
\begin{equation}
\label{fibers}
f^{-1}(z)= \{f^{-1}_1(z), f^{-1}_2(z), \dots , f^{-1}_d(z) \}
\end{equation}
where $d$ is the degree of $f$. For  $z\neq c_i, \infty$ define the (reduced) zero-homology group
$$
\tilde{H}_0(f^{-1}(z),\mathbb{Z})= \{n_1f^{-1}_1(z)+ n_2f^{-1}_1(z)+\dots n_df^{-1}_d(z): \sum n_i = 0, n_i \in \mathbb{Z} \}.
$$
It is a free $\mathbb{Z}$-module generated by 
$$f^{-1}_1(z)-f^{-1}_d(z), f^{-1}_2(z)-f^{-1}_d(z), \dots, f^{-1}_{d-1}(z)-f^{-1}_d(z)$$ and its dual space is denoted by $\tilde{H}^0(f^{-1}(z),\mathbb{C})$. The map (\ref{fmap}) induces homology and co-homology bundles with the base $\mathbb{C}\setminus \{c_1,\dots,c_k\}$ and fibers $\tilde{H}_0(f^{-1}(z),\mathbb{Z})$ and $\tilde{H}^0(f^{-1}(z),\mathbb{C})$. The continuous families of cycles
$$
f^{-1}_i(z)-f^{-1}_j(z)\in \tilde{H}_0(f^{-1}(z),\mathbb{Z})
$$
generate a basis of locally constant sections of a canonical connection on the homology bundle (the Gauss-Manin connection). Clearly a meromorphic function $g$ on $R$ defines a meromorphic section of the co-homology bundle, and  we may define a zero-dimensional Abelian integral as follows (see \cite{gamo07}).
 \bd
 \label{defai}
A zero-dimensional Abelian integral is an algebraic function
\begin{equation}
%\label{integral}
\int_{\delta(z)} g = n_1g(f^{-1}_1(z))+ n_2g(f^{-1}_1(z))+\dots + n_d g(f^{-1}_d(z)),
\end{equation}
where $g$ is a meromorphic function on $R$ and
\begin{equation}
\label{delta}
\delta(z) =  n_1f^{-1}_1(z)+ n_2f^{-1}_2(z)+\dots n_df^{-1}_d(z) \in \tilde{H}_0(f^{-1}(z),\mathbb{Z}) 
\end{equation}
is a continuous family of 0-cycles. 
\ed
Clearly, the functions $\phi_i$ in 
\e{generating} may be interpreted as zero-dimensional Abelian 
integrals
\be
%\label{integral}
\varphi_i(z) = \int_{\delta_i(z) }\frac{\omega}{df},
\ee
where 
\begin{equation}
\label{78}
\delta_i(z) = \sum_{j=1}^d n_{ij} f^{-1}_{j}(z)
\end{equation} 
and $n_{ij}$ are suitable integers (to be computed bellow).

To resume, we proved that if the set of poles of $\omega$ is contained in the set of poles of $f$, then the following statement is true. 

\begin{theorem}
\label{th1} The moments 
\be \l{mom}
m_s=\int_{\gamma} f^s \omega , \ \ \ s \geq 0
\ee
vanish
if and only if the zero-dimensional Abelian integrals
$$
\varphi_i(z) = \int_{\delta_i (z)}\frac{\omega}{df} , \ \ \ i= 1,2, \dots, k
$$
vanish.
\end{theorem}
Of course, to apply the Theorem \ref{th1} we need 
the precise values of the integer numbers $n_{ij}$ defined by (\ref{78}). Following \cite{mom}, we may compute these numbers as follows.
Consider
the pre-image of the star $S$ under $f$ 
$$
\lambda_f = f^{-1}(S) \subset R
$$
as a graph embedded in the Riemann surface $R$.
This graph, called a constellation, in a sense is a ``combinatorial portrait'' of the corresponding covering
(see \cite{lz} for details and different versions of this construction).
We will show that the combinatorial properties of $\lambda_f$ determine $n_{ij}$.

\begin{figure}[htbp]
%\medskip
%\epsfxsize=10.5truecm
%\centerline{\epsffile{1_1.eps}}
\begin{center}
\includegraphics[width=12cm]{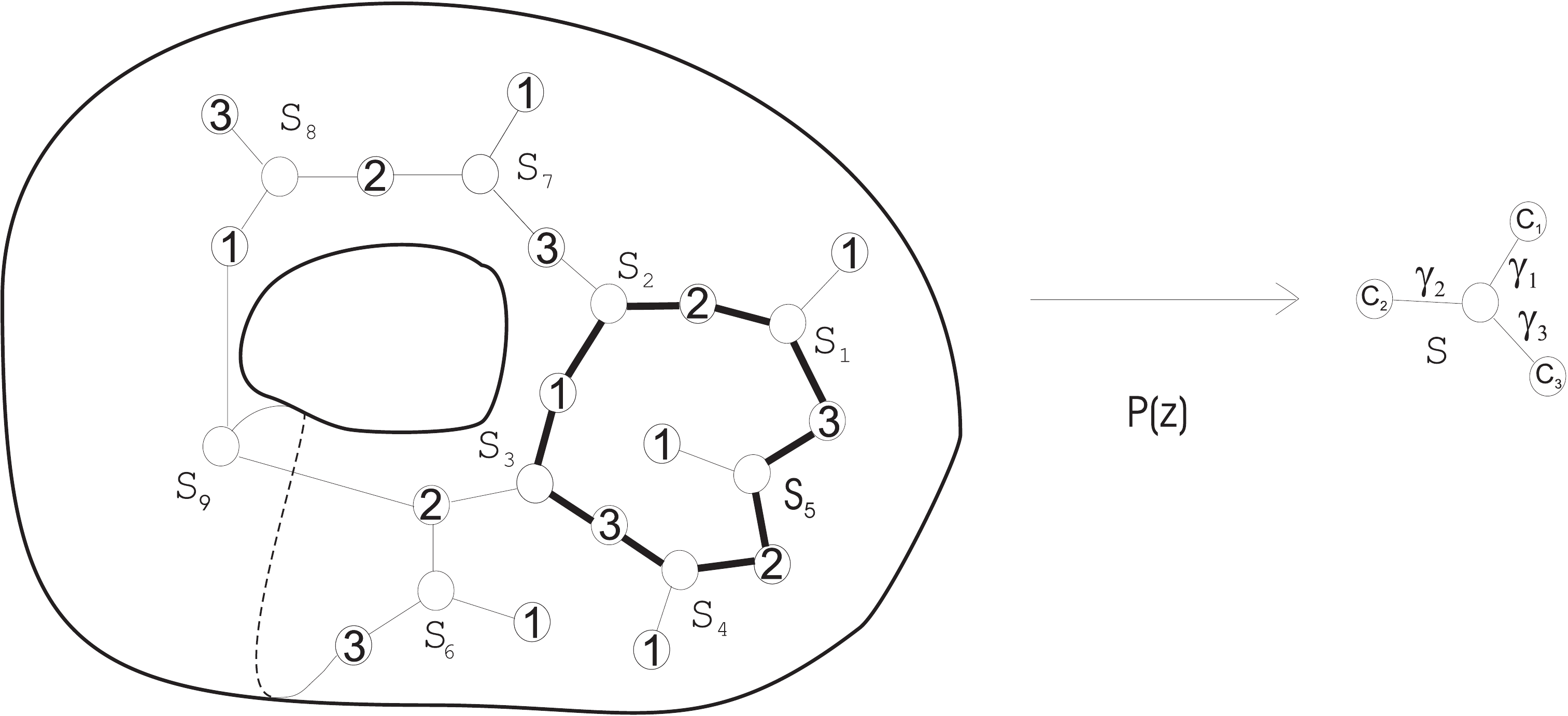}
\end{center}
\caption{}
\label{figure1}
%\smallskip
%\centerline{Figure 1}
%\medskip
\end{figure}

By construction, the restriction of $f(z)$ on $R\setminus \lambda_f$ is a covering of the topological punctured disk $\C\P^1\setminus \{S\cup \infty\}$ and therefore $R\setminus \lambda_f$ is a disjoint union of disks. This implies that the graph $\lambda_f$ is connected and the faces of $\lambda_f$ are in a one-to-one correspondence with poles of $f(z)$. 
For each $s,$ $1\leq s \leq k,$ we will mark vertices of $\lambda_f$ which are preimages of the point $c_s$ by the number $s$ (see Fig. 1).
Define a {\it star} of $\lambda_f$ as a subset of edges of $\lambda_f$ consisting of edges adjacent to some 
non-marked vertex. If $U$ is a simply-connected domain such 
that 
$S\setminus\{c_1,c_2,
... ,c_k\}\subset U$, then the set of stars
of $\lambda_P$ may be naturally identified with the set 
of single-valued branches of 
$f^{-1}(z)$ 
in $U$ as follows:
to the branch $f^{-1}_i(z),$ $1\leq i \leq n,$  
corresponds the star $S_i$
such that $f^{-1}_i(z)$ maps bijectively the interior of $S$ to the interior of $S_i$.

Now in order to obtain a deformation of  
the integration path in \e{generating0} 
satisfying the requirements above we may just to deform 
the initial path inside of each disk 
$R\setminus \lambda_f$ avoiding poles of $f$ (see Fig. 2).
\begin{figure}[ht]
%\medskip
%\epsfxsize=11.2truecm
%\centerline{\epsffile{1+_1.eps}}
%\smallskip
%\centerline{Figure 2}
\begin{center}
\includegraphics[width=12cm]{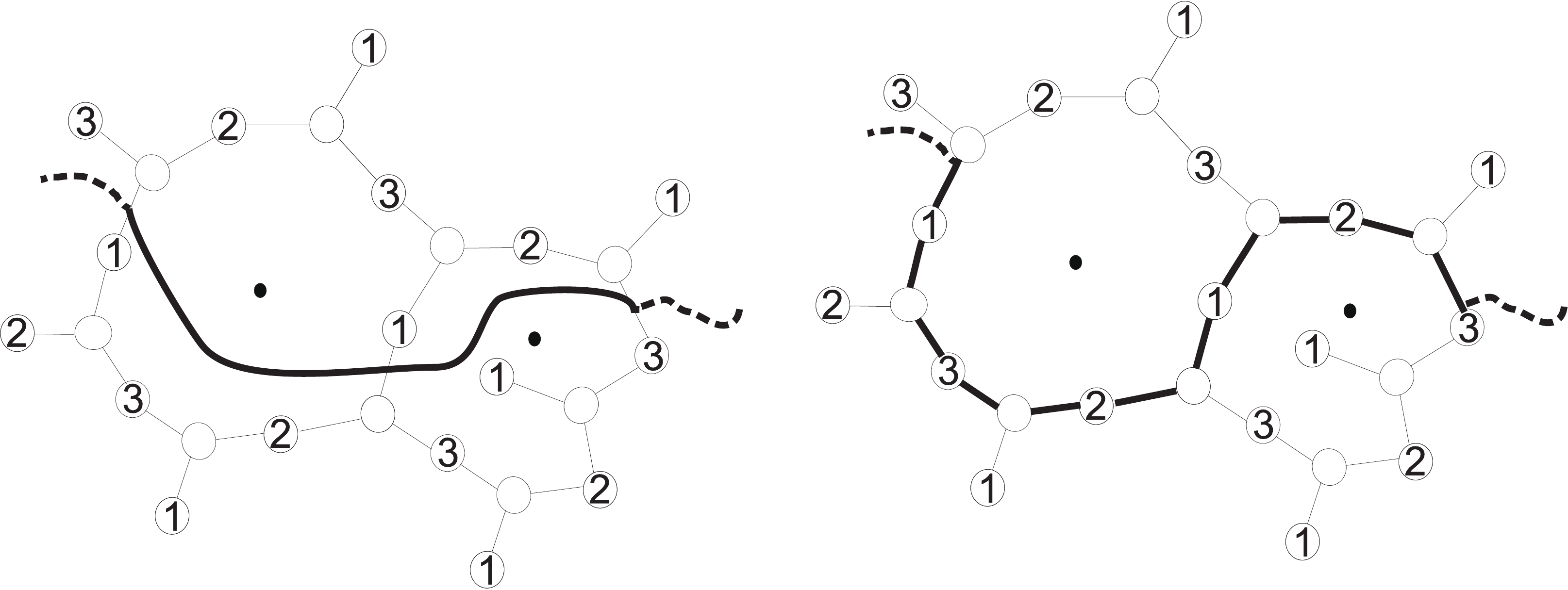}
\end{center}
\caption{}
\label{figure2}
%\medskip
\end{figure}
%\noindent 
If $\gamma$ is such a deformation and
$c_{s,i}$ is a unique vertex of the star $S_i$ marked by the number $s,$ then 
for $s,$ $1\leq s \leq k,$ we have:\be  
\l{piz}
\phi_s(z)=\sum_{i=1}^n n_{si}\left(\frac{\omega}{df}\right)(f^{-1}_{i}(z)), 
\ee
where $n_{si}$ is a sum of ``signed'' appearances of 
the vertex $c_{si}$ 
on the path $\gamma$.
By definition, this means that an appearance
is taken with the sign plus if the center of $S_i$ is followed by $c_{si}$, and minus
if $c_{si}$
is followed by the center of $S_i$.
For example, for the graph $\lambda_f$ shown on Fig. 1 and 
the path $\gamma\subset \lambda_f$ pictured by the fat 
line we have:  
\begin{gather*}
\phi_1(z)=\left(\frac{\omega}{df}\right)
(f^{-1}_{3}(z))-\left(\frac{\omega}{df}\right)(f^{-1}_{2}(z)),\\
\phi_2(z)=\left(\frac{\omega}{df}\right)(f^{-1}_{2}(z))-\left(\frac{\omega}{df}\right)(f^{-1}_{1}(z))+\left(\frac{\omega}{df}\right)(f^{-1}_{5}(z))- \left(\frac{\omega}{df}\right)(f^{-1}_{4}(z)),\\
\phi_3(z)=\left(\frac{\omega}{df}\right)(f^{-1}_{1}(z))-\left(\frac{\omega}{df}\right)(f^{-1}_{5}(z))+ \left(\frac{\omega}{df}\right)(f^{-1}_{4}(z))-\left(\frac{\omega}{df}\right)(f^{-1}_{3}(z)).
\end{gather*}

Clearly, since $\gamma$ is a closed loop, 
the above construction implies that $\sum_j n_{ij} =0$, so 
the cycles $\delta_i$ in \e{78} are contained in $\tilde{H}_0(f^{-1}(z),\mathbb{Z})$ indeed. Furthermore, 
since $R\setminus \lambda_f$ is a disjoint union of disks each of which contains a single pole of $f$, 
the following statement holds.

\bc \l{c1} If the curve $\gamma$ is not homological to zero in $R$ with poles of $f$ removed, then 
the vanishing of moments \e{mom} implies that there exists 
a non-zero cycle $\delta \in \tilde{H}_0(f^{-1}(z),\mathbb{Z})$
such that $\int_{\delta}\frac{\omega}{df}=0.$
\ec

Observe that Theorem \ref{th1} and Corollary \ref{c1} remain true 
without the restriction that the set of poles of $\omega$ is contained in the set of poles of $f$ if 
to change the condition $J(t)\equiv 0$ to the condition that 
$J(t)$ is rational.
Indeed, we always may find a polynomial $R$ 
such that the set of poles of the form $\widetilde \omega=R(f)\,\omega $ is contained in the set of poles of $f$. On the other hand, it follows from the definition of
$J(t)$ as a generating function that 
the functions $J(\omega,f,\gamma,t)$ and 
$J(f\omega,f,\gamma,t)$ are related 
by the equality 
$$J(f\omega,f,\gamma,t)=J(\omega,f,\gamma,t)t-\int_{\gamma}\omega$$
which implies inductively that the function 
$J(\omega,f,\gamma,t)$
is rational 
if and only if the function $J(\widetilde\omega,f,\gamma,t)$ does.

Further, observe that the above method may be applied also in the situation where the curve $\gamma$ is not closed and/or
is not connected 
(see \cite{mom}, Section 3
for the rational case which extends to the general case in the same way as above). Of coarse, if $\gamma$ is 
non-closed, then the condition $\sum_j n_{ij} =0$ for
$\delta_i$ in \e{78} is not necessary true.

\subsection{Solution of the moment problem in the case of generic position} 

Let $f$ be a meromorphic function on a compact Riemann surface  $R$, $z_0$ be a fixed regular value of $f$, and 
$\Delta$ be the set of critical values of $f$.
Recall that the monodromy group $G_f$ of the function $f$ is defined as 
the image of the homomorphism 
\begin{equation}
\label{representation}
\pi_1(\C \setminus \Delta, z_0) \rightarrow Aut(f^{-1}(z_0)),
\end{equation}
where $Aut(f^{-1}(z_0))$ is the full permutation group.
Further, a meromorphic function 
$f:\, R\rightarrow \C\P^1$ can be decomposed into a composition $f=p\circ q$ of holomorphic 
functions $q:\, R\rightarrow C$ and $p:\, C\rightarrow \C\P^1$, where $C$ is another compact Riemann surface, 
if and only if the group $G_{f}$ has an imprimitivity system consisting of 
$d=\deg p$ blocks. 
Observe that if the set $\{1,2,\dots, n\}$ is identified with the set $f^{-1}\{z_0\}$
then the set of blocks of the imprimitivity system corresponding to the decomposition  
$f=p\circ q$  
has the form 
$\f B_i=q^{-1}\{t_i\},$ $1\leq i \leq d_1,$ where $\{t_1,t_2,\dots, t_{d_1}\} =p^{-1}\{z_0\}$.
Finally, notice that  
if $f=\tilde p\circ \tilde q$, where  
$\tilde q:\, R\rightarrow \tilde C,$ $\tilde p:\, \tilde C\rightarrow \C\P^1$,
is an other decomposition of $f$ then  
the corresponding imprimitivity systems coincide if and only there exists an isomorphism $\mu:\, \tilde C \rightarrow C$
such that $$ p=\tilde p \circ \mu^{-1}, \ \ \ q=\mu \circ \tilde q.$$
In this case the decompositions $p\circ q$ and $\tilde p\circ \tilde q$
are called equivalent.

We say that meromorphic functions $f,$ $g$ on a Riemann surface $R$
have a non-trivial common compositional right factor
if there exists a Riemman surface $\tt R$, a holomorphic function $h\,:\,  R\rightarrow \tt R$ of degree greater than one, and holomorphic functions 
$\tt f, \tt g\,:\, \tt R\rightarrow \C\P^1$ such that $f=\tt f\circ h,$ $g=\tt g\circ h.$ The property of two
functions $f,$ $g$ to have a common compositional right factor
may be expressed via the vanishing of some zero-dimensional Abelian integrals. 
\bp \l{ecc}
Two meromorphic functions
 $f,$ $g$ on a compact Riemann surface $R$ have a common compositional right factor if and only if there exists a cycle $\delta(z) \in \tilde{H}_0(f^{-1}(z),\mathbb{Z})$
of the form $f^{-1}_i(z)-f^{-1}_j(z)$ such that 
\be  \l{eee} \int_{\delta(z)}g \equiv 0.\ee 
In particular, equality \e{eee} holds for all $\delta(z) \in \tilde{H}_0(f^{-1}(z),\mathbb{Z})$ if and only if  
there exists a rational function $\tt g$ such that $g=\tt g\circ f.$
\ep
\pr It is easy to see by the analytical continuation that, for a fixed index $i$, the indices $j$ satisfying the equality
\be \l{ee} g(f^{-1}_i(z))=g(f^{-1}_j(z))\ee 
form a block of an imprimitivity system $I$ with respect to the action 
of the monodromy group $G_f$ of $f$ on fibers of $f$. If \e{eee} holds, then this block contains more than one element and 
therefore there exists $\tt R$ and a meromorphic function $h\,:\, R\rightarrow \tt R$ such that $f=\tt f\circ h$ and fibers of $h$ coincide with 
blocks of $I.$  
Furthermore, equalities \e{ee}
imply that the function $g$ is constant on fibers of $h$. Therefore, a 
function $\tt g =g\circ h^{-1}$ is well defined and satisfies the equality 
$g=\tt g\circ h.$ \qed
\vskip 0.2cm

Notice that Proposition \ref{ecc} permits to reduce the moment 
problem for a collection $f$, $\omega,$ $\gamma$, $R$ to a similar problem for another collection $\tt f$, $\tt \omega,$ $\tt \gamma$, $\tt R$ with $\deg \tt f<\deg f$ 
whenever the space generated by cycles (\ref{78}), 
 contains a simple cycle  $f^{-1}_i(z)-f^{-1}_j(z)$, and in certain cases such a reduction may imply a complete solution of the initial problem
see e.g. \cite{let}, \cite{bull}, \cite{mom}, \cite{pry04}, \cite{krma10}.
In the particular case $R=\C\P^1$ the above Proposition is well known (see e.g. \cite{let}, Corollary 1) and follows easily from the L\"{u}roth theorem.
%\er
In particular, Proposition \ref{ecc} permits to solve the moment 
problem for $f$ in a generic position. 

\bt \l{yy} If the monodromy group $G_f$ of $f$ is the full symmetric group of $n=\deg f$ elements, then 
the vanishing of moments \e{mom} implies that 
either $\gamma$ is homological to zero in $R$ with poles of $f$ removed, or there 
exists a rational function $Q$ such that $\omega=Q(f)\,d f$ and 
$f(\gamma)$ is homological to zero in $\C\P^1$ with poles of $Q$ removed.
\et 

\pr
If $\gamma$ is not homological to zero in $R$ with poles of $f$ removed, then 
by Corollary \ref{c1} there exist integer numbers 
$\alpha_1,\alpha_2,\dots, \alpha_n$ 
such that 
\be \l{sub} \sum_{i=1}^n\alpha_{i}\left(\frac{\omega}{df}\right)(f^{-1}_{i}(z))=0\ee and $\sum_{i=1}^n\alpha_i=0.$ 
The last equality implies that the numbers $\alpha_1,\alpha_2,\dots, \alpha_n$ are not all  
equal between 
themselves. Let us assume that $\alpha_1\neq \alpha_2.$ 

Since  $G_f$ is a full symmetric group 
it contains the transposition $(12)$. Applying it to the equality \e{sub} and subtracting we obtain the equality $$(\alpha_1-\alpha_2)\left(\left(\frac{\omega}{df}\right)(f^{-1}_{1}(z))-
\left(\frac{\omega}{df}\right)(f^{-1}_{2}(z))\right)=0$$ implying the equality \be \l{g} \left(\frac{\omega}{df}\right)(f^{-1}_{1}(z))=\left(\frac{\omega}{df}\right)(f^{-1}_{2}(z)).\ee 
Since the full symmetric group of $f$ 
is primitive, the function $f$ is indecomposable. Therefore, 
Proposition \ref{ecc} applied to equality \e{g} implies that 
there exists a rational function $Q$ such that $\frac{\omega}{df}=Q(f).$
Hence, moments \e{zer} equal to the moments
$$\int_{f(\gamma)}Q(z)z^sdz, \ \ \ s\geq 0,$$ and the statement
follows from the classical result of the complex analysis. \qed

\vskip 0.2cm
Notice that Theorem \ref{yy}
also may be deduced from the characterization of doubly transitive groups
via the structure of their irreducible subspaces over $\C$, see \cite{let}, \cite{mom}
where this approach is used for rational $f$ and $g.$

\section{Vanishing of genus zero Abelian integrals} 
\label{abint}
In this section we give necessary and sufficient conditions
for an Abelian integral of dimension zero  to vanish identically in the case where $R$ is the Riemann sphere 
and the functions $f$ and $g$ are polynomials. More precisely, 
we solve the following problem:\\
\begin{quotation}
\emph{For a given polynomial $P$ of degree $n$ and 
a cycle $\delta(z)\in \tilde{H}_0(P^{-1}(z),\mathbb{Z})$
describe the polynomials $Q$ such that the associated Abelian integral }
\be \l{1} I(z)=\int_{\delta(z)}Q=\sum_{i=1}^n m_iQ(P^{-1}_i(z))\ee
\emph{vanishes identically.}\\
\end{quotation}
In distinction with the previous section we will not assume that 
$\tilde{H}_0(P^{-1}(z),\mathbb{Z})$ is reduced. 
Thus $\delta(z)$ may be any expression of the form
$$\delta(z)= v_1 P^{-1}_1(z)+v_2P^{-1}_2(z) + \dots + v_n P^{-1}_n(z),$$ where $v_i\in \Q.$ It is convenient to identify 
the cycle $\delta(z)$ with the vector 
$$\vec \delta=(v_1, v_2, \dots, v_n)$$ of $\Q^n.$ Under such an identification the natural action of $G_P$ on $\tilde{H}_0(P^{-1}(z),\mathbb{Z})$ descends to an action of $G_P$ of $\Q^n$ defining a
{\it permutation representation} of the group $G_P$  
\begin{equation}
\label{prepresentation}
\rho\, :\,   G_P \rightarrow GL(\Q^n).
\end{equation}

The understanding of irreducible components of $\rho$ plays a crucial role in the solution of the problem above. Indeed, let $Z_{\delta}$ be the vector space consisting of polynomials $Q$ such that Abelian integral \eqref{1} vanishes 
identically, and $V_\delta$ be the minimal $\rho$-invariant vector subspace 
of $\Q^n$ containing the vector $\vec \delta$. Then  
it is easy to see by the analytical continuation that $\int_{\delta(z)}Q\equiv 0$ if and only if $\int_{\gamma(z)}Q\equiv 0$ for any $\gamma (z)\in\tilde{H}_0(P^{-1}(z),\mathbb{Z})$ such that $\vec \gamma \in  V_\delta$.
This implies that in order to describe $Z_{\delta}$ it is enough 
to solve the following three problems:
\begin{enumerate}
\item First,  describe all possible irreducible $G_P$-invariant subspaces of $\Q^n$.
\item Second, for a given $\vec \delta \in \Q^n$  provide a method which allows to decompose the invariant space $V_{\delta}$ into a direct sum of 
irreducible $G_P$-invariant subspaces.

\item Third, for a given irreducible $G_P$-invariant subspace $V$,
describe the vector space $Z_{V}$ consisting of polynomials $Q$ 
such that $\int_{\delta(z)}Q\equiv 0$ 
for all $\delta \in V$.
\end{enumerate}

The solution of these problems is implicitly contained in \cite{pm} as an ingredient of the proof of the so called polynomial moment problem. For the reader's convenience we present to the end of the section  a self-contained proof, except of the classification irreducible $G_P$-invariant subspaces which is proved in \cite{pm} in a closed form.

\subsection{Description of the irreducible $G_P$-invariant subspaces of $\Q^n$. }

Below we reproduce the description of $G_P$-invariant subspaces of $\Q^n$ obtained in \cite{pm}.
More generally, we will describe 
$G$-invariant subspaces of $\Q^n$ for a permutation representation $G$ of an arbitrary
permutation group $G \subset S_n$ containing a cycle of length $n$
(the monodromy group of a polynomial of degree $n$ always contains such a cycle which corresponds to a loop around infinity). Notice that if $P$ is decomposed into a 
composition of rational functions $P=A\circ B$, then 
in the corresponding equivalence class there exists 
a decomposition into a composition of {\it polynomials} and below we always will consider only such decompositions.
In order to keep the correspondence between imprimitivity 
systems of $G_P$ and equivalence classes of decompositions 
of $P$ we need to modify
the definition of equivalent decompositions correspondingly.
Namely, we will call decompositions  
$P=A_1\circ W_1$ and $P=A_2\circ W_2$ equivalent if there exists a {\it linear function} $\nu$  such that $A_2=A_1\circ \nu,$ $W_2=\nu^{-1}\circ W_1.$ Abusing of notation, usually we will mean by a decomposition a corresponding equivalence class of decompositions.

Without loss of generality we may assume that a cycle of length $n$ contained  
in $G$ coincides with the cycle $(1,...,n)$. Then it is easy to see
that any imprimitivity system for $G$ coincides with the residue classes modulo $d$ for some $d\,|\,n$. In particular, this implies that any two decompositions $P=P_1\circ W_1$ and $P=P_2\circ W_2$ 
with $\deg P_1=\deg P_2$ are equivalent.
For each $d\,|\,n$ we denote by $V_d$ the subspace of $\Q^{n}$
consisting of $d$-periodic vectors of the form
$$
(v_1,\dots,v_d,v_1,\dots,v_d,\dots,v_1,\dots,v_d).
$$
It is easy to see 
that for given $d$ residue classes modulo $d$
form an imprimitivity system for $G$
if and only if the subspace $V_d$ is $G$-invariant.

Denote by $D(G)$ the set of all divisors of $n$ for which $V_d$ is $G$-invariant. In view of the above remarks $D(G)$ consists of numbers $d$ for which there exists a decomposition $P=A\circ W$ with $\deg A=d.$
Notice that $D(G)$ is a
lattice with respect to the ope\-rations $\land, \lor,$ where
$d\land \tilde d := \gcd(d,\tilde d)$ and $d\lor \tilde d := {\rm lcm}(d,\tilde d)$. Indeed, for an element $x\in X$
the intersection of two blocks containing $x$ and corresponding to $d,\tilde d\in D(G)$ is a block which
corresponds to $d\lor \tilde d$\,. On the other hand,
the intersection of two invariant subspaces $V_d,V_{\tilde d}$ is
an invariant subspace which is equal to $V_{d\land \tilde d}$\,.
We say that $d\in D(G)$ {\it covers} $\tilde d\in D(G)$
if $\tilde d \,|\,d,$ $\tilde d<d,$ and there exists no $l\in D(G)$ such that $\tilde d<l<d$ and $\tilde d\,\vert l,$ $l\vert d$. 

\bt \l{t1} {\rm (\cite{pm})} Assume a permutation group $G\subseteq S_n$ contains the cycle $(1\,2\,...\,n).$ Then 
each $G$-irreducible invariant subspace of $\Q^{n}$ has the form
\begin{equation}
\label{ud}
U_d:=V_d\cap \left(V_{d_1}^\perp\cap ...\cap V_{d_\ell}^\perp\right)
\end{equation}
where $d\in D(G)$ and
$d_1,...,d_\ell$ is a complete set of elements of $D(G)$ covered by $d$. The subspaces $U_d$ are mutually orthogonal and
every $G$-invariant subspace of $\Q^{n}$ is a direct sum of some $U_d$ as above. 
\et

\subsection{Decomposition of $V_{\delta}$ into a direct sum of irreducible subspaces}

The vectors
\begin{equation}
\label{wk}
\vec w_k=(1,\varepsilon_n^k,\varepsilon_n^{2k}, \,...\,,
\varepsilon_n^{(n-1)k}), \quad  \varepsilon = e^{2\pi i /n},\ \ \ 1\leq k \leq n
\end{equation}
form an orthonormal basis of $\C^n$ 
with respect to the standard Hermitian scalar vector product in $\C^n$,
and for any divisor $d$ of $n$ the vectors $\vec w_k$ for which $(n/d)\,\vert \,k$
form a basis of the complexification $V_d^{\C}$ of the subspace $V_d$. 
Furthermore, defining $\Psi_d$, $d\in D(G)$, as a subset of $\{1,2,\dots,n\}$ consisting of numbers $r$ such that $n/d$ is a divisor of $r$ but for any element $\tilde d \in D(G)$ covered by $d$ the number $n/\tilde d$
is not a divisor of $r$, we see that the vectors $\vec w_r,$ $r\in \Psi_d,$ form a basis 
of $U_d^{\C}.$ 
\bt \l{t2}  
The subspace 
$U_{d}$, $d\in D(G)$, is a component in the decomposition of the subspace $V_{\delta}$ 
into a sum of irreducible $G$-invariant subspaces of $\Q^n$, if and only if 
there exists a number $r\in \Psi_d$ such that $(\vec \delta, \vec w_{r})\neq 0.$
\et 
\noindent{\it Proof} (cf. Proposition 4.1 in \cite{pm}).
In view of Theorem \ref{t1}, it is enough to show that $V_{\delta}$ is orthogonal to $U_d$ if and only if $( \vec \delta, \vec w_{r})= 0$ for any 
$r\in \Psi_d$.
Let $V_{\delta}=\oplus\ U_{d_j}$ be the decomposition of $V_{\delta}$ 
into a sum of irreducible $G$-invariant subspaces. 
If $V_{\delta}$ is orthogonal to $U_d$, then
 $V_{\delta}$ should be orthogonal also to $U_d^{\C}$ implying that $(\vec \delta, \vec w_{r})= 0$ for all $r\in \Psi_d$. On the other hand, if 
$(\vec \delta, \vec w_{r})= 0$ for all $r\in \Psi_d$, then
$\vec \delta$ is orthogonal to $U_d$. Therefore, 
 $V_{\delta}$ is orthogonal to $U_d$ in view of the minimality of $V_{\delta}.$\qed

\subsection{Description of the spaces $Z_{U_d}$.}
The structure of $Z_{U_d}$ or, more generally, of any $Z_{\delta}$ 
is closely related to the compositional algebra of polynomials. 
For example, if 
$\delta(z) \in H_0(P^{-1}(z),\Q)$, then for any polynomial $Q$ of the form 
$ Q=R \circ P$ we have 
\be \l{ccoo} \int_{\delta(z)} Q=R(z) \sum_{i=1}^n v_i\ee 
implying that the corresponding integral vanishes whenever 
$\delta(z)$ is contained 
in the reduced homology group, or equivalently the vector $\vec \delta$ is contained in $V_1^{\perp}.$ 

Further, if $P=A\circ W,$ $\deg A=d,$ is a decomposition of $P$ corresponding to $d\in D(G_P)$, 
then 
for any branch $P_i^{-1}(z)$ of $P^{-1}(z)$ there exist a branch $W_j^{-1}(z)$ of $W^{-1}(z)$ and a branch $A_k^{-1}(z)$ of $A^{-1}(z)$ such that
\begin{equation}
\label{co}
P_i^{-1} = W_j^{-1} \circ A_k^{-1}.
\end{equation}
Therefore, for any 
cycle $\delta(z)\in H_0(P^{-1}(z),\Q)$ and polynomial $Q$ we have
\be 
\l{cco} 
\int_{\delta(z)}Q=\sum_{k=1}^d\left(\int_{\delta_{k,W}(z)}Q\right)\circ A_{k}^{-1}, 
\ee where $\delta_{k,W}(z)\in H_0(W^{-1}(z),\Q)$, $d=\deg A,$ and hence the integral $\int_{\delta(z)}Q$ 
vanishes identically whenever all the integrals $\int_{\delta_{k,W}(z)}Q$, $1\leq k \leq d,$ do. In view of the above remark, this implies in particular
that if $\delta(z)\in H_0(P^{-1}(z),\Q)$ is a cycle such that all cycles $\delta_{k,W}(z)$, $1\leq k\leq d,$ are in the reduced homology group
$\tilde{H}_0(W^{-1}(z),\Q)$, then $\int_{\delta(z)}Q$ vanishes identically for any polynomial $Q$ of the form  
$Q=B\circ W$.

In the following we always will assume that the numeration of roots $P^{-1}_i(z)$ of $P^{-1}(z)$ satisfies the requirement that the cycle in $G_P$ corresponding to a loop around 
infinity in $\C$ coincides with the cycle
$(1\,2\,\dots\, n).$ In particular, such a choice of the numeration yields that 
when $i$ in formula \e{co} runs a residues class modulo $d$, corresponding $j$ runs numbers
from $1$ to $n/d$ while $k$ remains fixed, implying that 
%the vectors $\vec \delta_{k,W}\in \Q^{n/d}$ correspodning to 
the cycles $\delta_{k,W}(z)\in H_0(W^{-1}(z),\Q)$ in \eqref{cco} are reduced if and only for and $k,$ $1\leq k \leq d,$ 
the equality  
\be \l{wert} (\vec \delta, \vec e_{k,d})=0\ee holds, 
where $\vec e_{k,d},$ $1\leq k \leq d,$ denotes a vector of $\Q^n$ with coordinates $v_1,v_2,\dots, v_n$ such that 
$v_i=1$ if $i=k\,\mod d$, and $v_i=0$ otherwise. Since vectors
$\vec e_{j,d},$ $1\leq j \leq d,$ obviously form a basis of 
$V_d$, we conclude that the cycles $\delta_{k,W}(z)\in H_0(W^{-1}(z),\Q)$ are reduced if and only if $\vec\delta$ is orthogonal to $V_d.$  

Returning to the description of the space $Z_{U_d}$ observe that it always contains the space $Z_{V_d}$ in view of the inclusion  
$U_d\subset V_d$. Furthermore, if $\tilde d$ if an element of $D(G_P)$ covered by $d$ and $P=\tilde A \circ \tilde W$ if a decomposition
corresponding to $\tilde d$, 
%and $\delta_{k,\tilde W}(z)$, $1\leq k \leq \tilde  d,$ be cycles 
%of $\tilde{H}_0(\tilde W^{-1}(z),\Q)$ corresponidng to the decomposition $P=\tilde A \circ \tilde W$. 
then, since 
$U_d$ is orthogonal to the subspace $V_{\tilde d}$, the cycles $\delta_{k,\tilde W}(z)$, $1\leq k \leq \tilde  d,$ are reduced. Therefore, for any such $\tilde d$, 
the ring $\C[\tilde W]$ of polynomials in $\tilde W$ is contained in the space $Z_{U_D}$.

 \bt \l{t3} Let $d$ be an element of $D(G_P)$. Furthermore, let 
$d_1,...,d_\ell$ be a complete set of elements of $D(G)$ covered by $d$ and 
$P=A_{i}\circ  W_{i},$ $1\leq i \leq \ell,$ be the corresponding decompositions. Then 
\be\l{es} Z_{U_d}=Z_{V_d}+\C[W_{1}]+\C[W_{2}]+\dots +\C[W_{l}].\ee 
\et 
\noindent{\it Proof} (cf. Theorem 1.1 in \cite{pm}).
In view of the above remarks, the right part of 
\e{es} is contained in $Z_{U_d}.$ So, we only must establish the inverse inclusion.

First observe that the numeration of branches of $P^{-1}(z)$ implies that at points close enough to infinity 
the functions $Q(P^{-1}_i(z))$, $1\leq i \leq n,$ 
may be represented by a converging Pusieux series \be \l{ps2}
Q(P^{-1}_i(z))=\sum_{k=-q}^{\infty}
s_k\varepsilon_n^{(i-1)k}z^{-\frac{k}{n}},
\ee where $q=\deg Q(z)$ and $\v_n=\exp(2\pi i/n).$ 
Furthermore, substituting \e{ps2} to \e{1} we see that 
the integral $\int_{\delta(z)}Q$ 
vanishes identically if and only if for any $k\geq -q$ the equality 
$$
\sum_{i=1}^{n}v_i s_k\varepsilon^{(i-1)k}_n=(\overline{\vec \delta},\vec w_{k})s_k=0
$$ holds. In particular, if $Q(z)\in Z_{U_d},$ $d\in D(G_P),$ then equalities 
\be \l{vot}
(\overline{\vec v},\vec w_{k})s_k=0, \ \ \ k\geq -q,
\ee 
hold for any $\vec v\in U_d$. Moreover, they hold for any $\vec v\in U_d^{\C}$. 
Since $U_d^{\C}$ is generated by the set of vectors $\vec w_r,$ $r\in \Psi_d$, and
this set transforms to itself under the complex conjugation, this implies that if $Q(z)\in Z_{U_d},$ then for any $r\in \Psi_d$ the equality
$s_k=0$ 
holds whenever $k\equiv r\, \mod n.$ Furthermore, 
clearly
the inverse is also true. 
Similarly, $Q(z)\in Z_{V_d}$ if and only if $s_k=0$ for any $k$ such that $(n/d)\vert k$.

Assume now that $Q(z)\in Z_{U_d}$ and consider series \e{ps2}.
If $s_k=0$ for any $k$ such that $(n/d)\vert k$, then $Q(z)\in Z_{V_d}.$ 
So, suppose that there exists 
$k$ such that $(n/d)\vert k$ but $s_k\neq 0$. Then $(n/\tilde d)\vert k$ for some $\tilde d$ covered by $d$ and we may assume that $\tilde d=d_1$.
Set 
$$\psi(z)=\sum_{\substack {k\geq -q \\ k\equiv 0\, \mod n/ d_1}}s_{k} z^{-\frac{k}{n}},$$
where $s_k,$ $k\geq -q,$ are coefficients of series \eqref{ps2}. Clearly, we have:
\be \l{iio} \left(\frac{n}{d_1}\right)\psi(z)=
Q(P^{-1}_{1}(z))+
Q(P^{-1}_{d_1+1}(z))+Q(P^{-1}_{2d_1+1}(z))
+ ... + Q(P^{-1}_{n-d_1+1}(z)).
\ee
Since the set of indices appearing in the right part of \e{iio} is a block, 
the function $\psi(z)$ is invariant with respect to the subgroup of $G_P$
which stabilizes $P^{-1}_{1}(z)$. Therefore, by the main theorem of Galois theory, $\psi(z)$ is contained in the field $\C(z)(P^{-1}_{1}(z))=\C(P^{-1}_{1}(z))$ and hence 
$\psi(z)=R_1(P^{-1}_1(z))$ for some rational function $R_1$.
Moreover, $R_1$ is actually a polynomial since the right part of \e{iio} may have a pole 
only at infinity. Further, since \e{iio} implies 
that $$R_1(P^{-1}_1(z))=R_1(P^{-1}_{d_1+1}(z))=R_1(P^{-1}_{2d_1+1}(z))=\dots =R(P^{-1}_{n-d_1+1}(z)),$$ it follows from \e{co} that 
$R_1$ is constant on fibers of $W_{1}$. Reasoning now as in  
Proposition \ref{ecc} we conclude that $R_1=S_1\circ W_{1}$ for 
some polynomial $S_1$ (cf. Lemma 4.3 in \cite{pm}). 

Define a polynomial $T_1(z)$ by the equality $$T_1(z)=Q(z)-R_1(z).$$
Then by
construction the Puiseux series of $T_1(P^{-1}_1(z))$
contains no non-zero coefficients with indices which are multiple of
$n/ d_1$. If $T_1(z)$ is contained in $Z_W$, then $$Q(z)=T_1(z)+S_1( W_{1}(z))$$ and we are done. Otherwise 
arguing as above we may find polynomials $R_2,$ $S_2$ such that $R_2=S_2\circ W_{2}$
and the Puiseux expansion of $T_2(P^{-1}_1(z)),$ where $$T_2(z)=T_1(z)-S_2(z),$$
contains no non-zero coefficients whose indices are multiple of
$n/ d_1$ or $n/ d_2$. It is clear that continuing this process we eventually will arrive to some $T_s(z)$ which is contained in $Z_W$ and therefore to a representation 
$$Q(z)=T_s(z)+S_1(W_1(z))+S_2(W_{2}(z))+\dots + S_l(W_{l}(z)).\ \ \ \Box$$

\vskip 0.2cm

In view of Theorem \ref{t3} in order to complete the description of the space $Z_{U_d}$ we only must
describe the space $Z_{V_d}$. 
Observe that the vectors
$\vec e_{j,d},$ $1\leq j \leq d,$ defined above satisfies the equality
$$(\vec e_{j,d}, \vec w_k)=\varepsilon_n^{k(j-1)}
(\vec e_{1,d}, \vec w_k).$$ Therefore, in order to check 
that \e{vot} holds for any $\vec v\in Z_{V_d}$ it is enough 
to check that it holds 
for one single vector $\vec e_{1,d}$. In other words, 
the space $Z_{V_d}$ consists of 
polynomials $Q(z)$ satisfying the equality 
\be \l{con}
Q(P^{-1}_{1}(z))+
Q(P^{-1}_{d+1}(z))+Q(P^{-1}_{2d+1}(z))
+ ... + Q(P^{-1}_{n-d+1}(z))\equiv 0.
\ee 
Furthermore, if $P=A\circ W$ is a decomposition corresponding to 
$d\in D(G_P)$, then in view of \eqref{co} equality \eqref{con}   
reduces to the equality 
 \be \l{con1}
Q(W^{-1}_{1}(z))+
Q(W^{-1}_{2}(z))+Q(W^{-1}_{3}(z))
+ ... + Q(W^{-1}_{n/d}(z))\equiv 0.
\ee 

The Newton formulae imply that whenever $\deg Q<\deg W$ the sum in the left hand side of \eqref{con1} is a constant. Therefore, if
$\mu_i=W^{-1}_i(c)$, $i=1,2,\dots,n/d$, for some generic $c$, then  the relation 
$$Q(\mu_1)+Q(\mu_2)+\dots +Q(\mu_{n/d})=0$$
describes the vector space of polynomials satisfying (\ref{con1}).
Furthermore, for $Q(z)$ of arbitrary degree it is easy to see, using $W$-adic decomposition, that the sum in \e{con1} is a polynomial, and that $Q(z)$ satisfies \e{con1} if 
and only if all coefficients in its $W$-adic decomposition satisfy it.

\subsection{Corollaries}
In this subsection we discuss some particular cases of the above results
which may be useful for applications.
\bp \l{p1} Let $P$ be an indecomposable polynomial. 
If an Abelian integral $\int_{\delta(z)}Q $,
where $Q$ is a non-zero polynomial and $\delta(z)\in H_0(P^{-1}(z),\Q)$, vanishes identically, then 
either $Q$ is a polynomial in $P$ and the cycle $\delta(z)$ is reduced, or there exists a rational number $a$ such that 
\be \l{cyc} \delta(z) = a(P_1^{-1}(z)+ P_2^{-1}(z)+ \dots + P_n^{-1}(z)).\ee
\ep
\pr It follows from the assumption $Q(z)\not\equiv 0$ that the subspace $V_{\delta}$
does not coincide with whole $\Q^n.$ Therefore,
by Theorem \ref{t1} either $V_{\delta}=U_1$ or $V_{\delta}=U_n$.
Obviously in the first case \eqref{cyc} holds while in the second $\delta(z)$ is reduced. Furthermore, by Theorem \ref{t3} in the second case $Q$ is contained in $Z_{U_n}=Z_{V_n}+\C[P].$ However, since $V_n=\mathbb Q^n$, the space $Z_{V_n}$ is trivial and therefore $Q$ is a polynomial in $P.$ 
\vskip 0.2cm
Notice that the conclusion of Proposition \ref{p1} holds for any polynomial $P(z)$ in generic position since decomposable polynomials obviously 
form a proper algebraic subset in the set of all polynomials of degree $n.$
Notice also that in order to prove 
Proposition \ref{p1} one can use instead of Theorem \ref{t1} a classical result, relating doubly transitivity of a group
with the structure of its permutation representation over $\C$, and Schur theorem, relating doubly transitivity and primitivity for group containing a transitive cyclic subgroup (see \cite{let} for such an approach).  
\vskip 0.2cm
The conclusion similar to the one in Proposition \ref{p1} is
true for arbitrary $P$ if to impose some limitations on $\delta(z)$.
\bp \l{p2} If an Abelian integral $\int_{\delta(z)}Q $,
where $Q$ is a non-zero polynomial and $\delta(z)\in H_0(P^{-1}(z),\Q)$, vanishes identically, and for any $d\in D(G_P),$ $d\neq 1,$ there exists $r\in \Psi_d$ such that $(\vec \delta,\vec w_r)\neq 0$, then $Q$ is a polynomial in $P$ and the cycle $\delta(z)$ is reduced.
\ep
\pr 
It follows from Theorem \ref{t1} and Theorem \ref{t2} that 
$V_{\delta}=U_1^{\perp}.$ In particular, the cycle $\delta(z)$ is reduced. Furthermore, $V_{\delta}$ contains
vectors $e_i-e_j,$ $1\leq i,j \leq n,$ implying that 
\be \l{rav} Q(P^{-1}_1(z))=Q(P^{-1}_2(z))=\dots=Q(P^{-1}_n(z)).\ee 
Now the statement follows from Proposition \ref{ecc}. \qed

\vskip 0.2cm
A weaker version of Proposition \ref{p2} can be formulated as follows
 \bp \l{p3} If an Abelian integral $\int_{\delta(z)}Q $,
where $Q$ is a non-zero polynomial and $\delta(z)\in H_0(P^{-1}(z),\Q)$, vanishes identically, and 
there exists $r\in \Psi_n$ such that $(\vec \delta, w_r)\neq 0,$
then $Q(z)$ may be represented in the form
\be \l{ec} Q(z)=S_1(W_1(z))+S_2(W_2(z))+\dots +S_l(W_l(z)),\ee where 
$S_1,S_2,\dots,S_l$ are polynomials and
$W_1,W_2,\dots W_l$ are polynomial compositional right factors of $P(z)$.
\ep
\pr It follows from Theorem \ref{t2} that $V_{\delta}$ contains $U_n$. Therefore, $Z_{\delta}\subseteq Z_{U_n}$ and the statement follows from Theorem \ref{t3}.  \qed

\vskip 0.2cm

Notice that in general it is not true that for any 
polynomial $W_i$ appearing in representation \e{ec} the corresponding cycles 
\be \l{pro} \delta_{k,W_i}(z)\in \tilde{H}_0(W^{-1}_i(z),\Q)\ee
are reduced.
On the other hand, if it is true, then the theorem is equivalent to the statement that the integral $\int_{\delta(z)}Q$ vanishes if and only if it is 
a sum of pull-backs of integrals along cycles homological to zero. Therefore, it is of interest to give conditions providing
that there exists a representation \e{ec} such that all cycles \eqref{pro} are reduced.

\bp\l{p4} 
Assume that an Abelian integral $\int_{\delta(z)}Q $,
where $Q$ is a non-zero polynomial and $\delta(z)\in H_0(P^{-1}(z),\Q)$ vanishes identically. Furthermore, assume that 
there exist $d_1,d_2, \dots d_l\in D(G_P)$ such that $(\vec \delta, w_r)= 0$ if and only if $r$ is a divisor of one of 
$d_1,d_2, \dots d_l$.
Then $Q(z)$ may be represented in the form \eqref{ec} where all cycles \eqref{pro} 
are reduced.
\ep 
\pr Since, the condition of the theorem implies that $V_{\delta}$ coincides with the orthogonal complement to the sum of  
$V_{n/d_1},V_{n/d_2}, \dots, V_{n/d_l}$ in $\mathbb Q^n,$ the proof is obtained by an obvious modification of the proof of Theorem \ref{t3}.
\qed

\vskip 0.2cm
\noindent{\it Remark.}\ The results similar to Propositions \ref{p3}, \ref{p4} (without the solution of the general problem) were obtained in the recent preprint 
\cite{stm11} (Thereom 2.2) where they also were deduced from Theorem \ref{t1} by the method of \cite{pm}. 
However, in \cite{stm11} they seem to be formulated not in their full generality. For example, the second part of Theorem 2.2 which is an analog 
of our Proposition \ref{p3} is formulated under an unnecessary condition that $r\in \Psi_n$ for which $(\vec w, \vec w_r)\neq 0$
is coprime with $n.$ 

\subsection{Polynomial moment problem on a system of intervals with weights}
Recall that the polynomial moment problem, recently solved in \cite{pm}, asks about conditions implying that a polynomial $Q$ satisfies the system of equations
\be \l{lsko} \int_a^b P^idQ=0, \ \ \ i\geq 0,\ee where $P$ is a given polynomial and $a,b$ are fixed complex numbers. It was shown in \cite{pm} that any solution $Q$ of this problem has the form \eqref{ec} where all $W_i$ satisfy the condition $W_i(a)=W_i(b)$. In the above notation the approach of \cite{pm} may be summarized as follows. First,
to construct, starting from the collection $P$, $a,$ $b$,  
a collection of cycles $\delta_j(z),$ $1\leq j \leq s,$ in $ 
\tilde{H}_0(P^{-1}_i(z),\Q)$ such that \e{lsko} holds if and only if 
\be \l{cde} \int_{\delta_i(z)}Q=0, \ \ \ 1\leq j \leq s.\ee Then to 
show that the corresponding invariant subspace generated by the cycles $\delta_j(z),$ $1\leq j \leq s,$ always contains a vector $\vec v$ such that
$(\vec v, w_r)\neq 0$ for some $r$ coprime with $n$ (this was done in \cite{bull} by means the so-called ``monodromy lemma'') and, finally, to use Proposition \ref{p3} 
repeatedly.

It is not hard to see that the results of the current section may be interpreted as a solution of the polynomial moment problem ``on a system of intervals with weights''.
More precisely, for any collection consisting of a polynomial $P$, 
complex numbers $a_i,b_i$, $1\leq i \leq l,$ and rational numbers $c_i,$ $1\leq i \leq l,$  
using approach of 
Section 1 (see also \cite{pm}, \cite{mom} where more attention to non-closed curves is given) one can construct 
a finite collection of cycles $\delta_j(z),$ $1\leq j \leq s,$ in $ 
\tilde{H}_0(P^{-1}_i(z),\Q)$ such that the equalities 
\be \l{mpr} c_1\int_{a_1}^{b_1}P^idQ+ c_2\int_{a_2}^{b_2}P^idQ+\dots+c_l\int_{a_l}^{b_l}P^idQ=0, \ \ \ i\geq 0,\ee hold if and only if equalities \e{cde} hold. 
%\be \l{cde} \int_{\delta_i(z)}Q=0, \ \ \ 1\leq j \leq s.\ee
%Of coarse, as above, for a solution of such a moment problem of importancy is not the system  $\delta_j(z),$ $1\leq j \leq s,$ itself but the corresponding 
%minimal
Since the results of this section provide a description of $Q$ satisfying \e{cde}, we obtain therefore a description of solutions of \e{mpr}. 
Notice however that since $G_P$-invariant subspace $V$ of $\Q^n$ generated by the vectors
$\vec \delta_j,$ $1\leq j \leq s,$ obtained from \e{mpr} may be arbitrary we can not expect to have 
such a good description of solutions as for problem \e{lsko}.

As a simple illustration, let us take $P$ equal to $T_6,$ where $T_n$ denotes the $n$th Chebyshev polynomial,
$T_n(cos\phi)=cos(n\phi).$ It is easily seen that $T_n$ has only two finite critical values and that 
the corresponding constellation has the form shown on the Fig. 3,
%\begin{figure}[ht]
%\medskip
%\epsfxsize=6.2truecm
%\centerline{\epsffile{zz.eps}}
%\smallskip
%\centerline{Figure 3}
%\medskip
%\end{figure}

\begin{figure}[htbp]
%\medskip
%\epsfxsize=10.5truecm
%\centerline{\epsffile{1_1.eps}}
\begin{center}
\includegraphics[width=8cm]{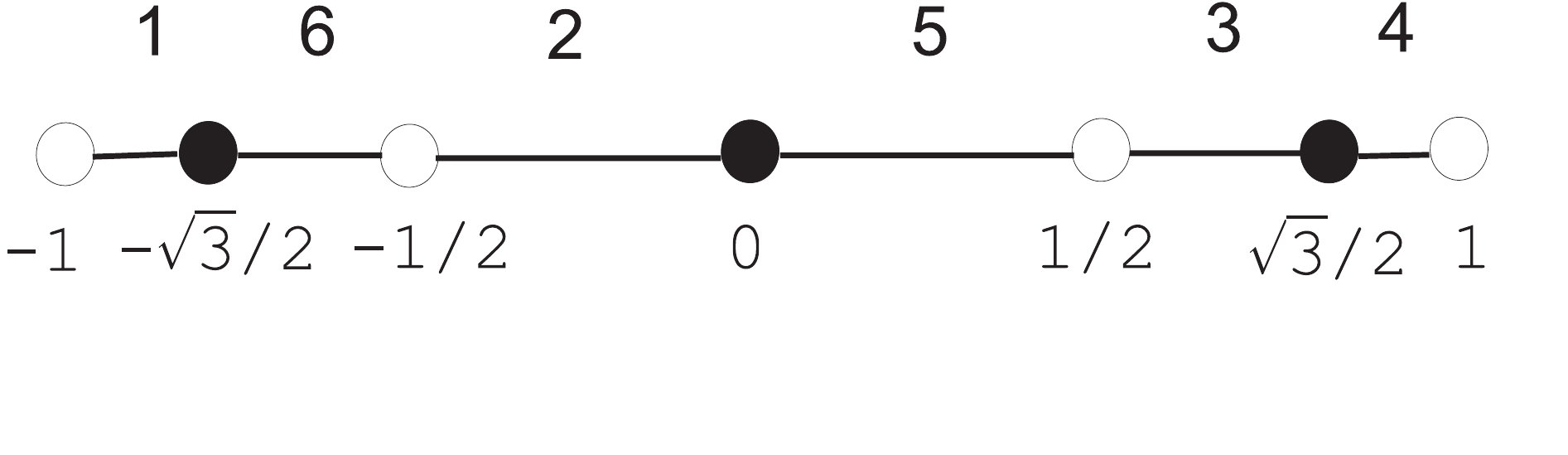}
\end{center}
\caption{}
\label{figure3}
%\smallskip
%\centerline{Figure 1}
%\medskip
\end{figure}
\noindent 
where the ``middle'' vertices of stars are omitted 
and the numeration of stars is
chosen in such a way that a permutation at infinity  
coincides with the cycle $(123456)$ (see e.g. \cite{ch}).

If we are searching for solutions of the moment problem 
\be \l{mprr} \int_{-\sqrt{3}/2}^{\sqrt{3}/2}T_6^idQ=0, \ \ \ i\geq 0,\ee 
on a single segment $[-\sqrt{3}/2, \sqrt{3}/2],$ 
we arrive to the vanishing problem for the Abelian integral 
$$\delta(z)=Q(T_{6,6}^{-1}(z))-Q(T_{6,2}^{-1}(z))+ 
Q(T_{6,5}^{-1}(z))-Q(T_{6,3}^{-1}(z)),$$ where by $T_{6,i}^{-1}(z),$ $1\leq i \leq 6,$ are denoted the branches of $T_{6}^{-1}$ (although we should obtain a cycle
for each critical value, it is easy to see that, since we have only two critical values,
the corresponding cycles are proportional).
Clearly, the corresponding 
vector $(0,-1,-1,0,1,1)\in \Q^6$ is orthogonal to both $V_2$ and $V_3$ implying that $V_{\delta}=U_6.$ Therefore, by Theorem \ref{t3} any solution of 
\e{mprr} has the form 
$$Q(z)=A(T_3(z))+B(T_2(z)),$$ where $A,B\in \C[z].$ 

On the other hand, the ``generalized'' moment problem 
\be \l{mprrr} \int_{-1}^{-1/2}T_6^idQ-\int_{-1/2}^{1/2}T_6^idQ+
\int_{1/2}^{1}T_6^idQ
=0, \ \ \ i\geq 0,\ee 
leads to 
the vanishing problem for the Abelian integral 
$$\delta(z)=Q(T_{6,1}^{-1}(z))-Q(T_{6,2}^{-1}(z))+Q(T_{6,3}^{-1}(z))
-Q(T_{6,4}^{-1}(z))+Q(T_{6,5}^{-1}(z))
-Q(T_{6,6}^{-1}(z)).$$
One can check that in this case the corresponding invariant subspace $V$ of $\Q^6$ 
coincides with $U_3$ and Theorem \ref{t3} implies that any solution of 
\e{mprrr} has the form 
$$Q(z)=A(T_6(z))+B(z),$$ where $A$ is an arbitrary polynomial and $B$ is a polynomial 
such that $$B(T_{3,1}^{-1}(z))+B(T_{3,2}^{-1}(z))
+B(T_{3,3}^{-1}(z))\equiv 0.$$

\section{Vanishing of hyperelliptic Abelian integrals}
Let $f\in \C[x]$ be a polynomial and $\Gamma_t= \{(x,y)\in \C^2: y^2-f(x)=t\}$ a family of hyperelliptic curves. Consider the Abelian integral
\begin{equation}
\label{abelian}
I(t)=\int_{\gamma(t)} \omega
\end{equation}
where $\omega = P(x,y)dx + Q(x,y) dy$ is a polynomial one form, and $\gamma(t) \in H_1(\Gamma_t,\mathbb{Z})$ is a continuous family of 1-cycles. 

The purpose of this section is to determine necessary and sufficient conditions for
the Abelian integral $I$ to be single valued, polynomial, or rational function. These three conditions are in fact equivalent. Indeed, $I(t)$ is a function of moderate growth, with a bounded modulus in any sector, centered at a singularity. Thus $I(t)$ is single-valued if and only if it is a rational, in fact polynomial function.
\subsection{Vanishing and the moment problem.}
The derivatives of $I$ can be seen as moments on a Riemann surface and this permits to apply the results of the preceding section. Indeed, every polynomial one-form $\omega$
can be written as
$$
 \omega= k(x)y dx + dA+ B d(y^2-f(x)), \quad A, B \in \C[x,y], k \in  \C[x] .
 $$
Therefore,
\begin{equation}
\label{abel}
I(t)=  \int_{\gamma(t)} k(x)y dx , \quad I'(t)= \frac{1}{2} \int_{\gamma(t)} \frac{k(x) dx }{y} 
\end{equation}
and more generally
\begin{equation}
\label{ik1}
I^{(k+1)}(t) = (1/2)(-1/2)(-3/2)\dots(-k+1/2)  \int_{\gamma(t)} \frac{k(x)}{y^{2k+1} }dx, k\geq 0 .
\end{equation}
Thus, 
$$
I^{(k+1)}(0)= m_k= (1/2)(-1/2)(-3/2)\dots(-k+1/2)  \int_{\gamma(0)} g^{k+1} \omega$$
 where  
 $$
 \omega=
  k(x) y dx,\ \ g = \frac{1}{f} ,
$$
implying that the Abelian integral $I(t)$ vanishes identically if and only if the moments $\int_{\gamma(0)} g^k \omega$, $k\geq 0$, vanish.
Furthermore, if to replace $\omega$ by $g^k\omega$, then, for $k$ sufficiently big, the set of poles of $g^k\omega$ will be a subset of the set of poles of $g$ and the results 
of Section 1 apply.

The zero-dimensional integrals described in Theorem \ref{th1} take the form
$$
\varphi_i(z) = \int_{\delta_i }\frac{\omega}{f^kdf} = \int_{\delta_i (z)}\frac{k(x)}{f^{(2k-1)/2}f'(x)}= z^{-(2k-1)/2}\frac{d}{dz} \int_{\delta_i (z)} K(x)
$$
where $K(x)= \int k(x)dx$ is a primitive of $k$, and the zero-cycles $\delta_i$ are constructed from the constellation $\lambda_f=f^{-1}(S)$
%, where $c_i$, $i\geq 0$ are the finite critical values of $1/f$ (and $0$ is a finite critical value corresponding to $x=\infty$), 
as explained in section \ref{section2}. The above gives necessary and sufficient conditions for the moments $m_i$, $i\geq k$ to vanish or, equivalently, for $I(t)$ to be a polynomial.
%As $I(t)$ can not have poles on the finite plane
Thus, we have proved
\begin{theorem}
\label{main3}
The Abelian integral  (\ref{abelian}) is a rational function if and only if the 0-dimensional integrals
$$\int_{\delta_i (z)} K(x), i=1,2, \dots, k$$ are identically constant.
\end{theorem}
\noindent{\bf Remark.}
Consider 
the polynomial $f(x)= (x^2/2-1)^2$ as on
fig.\ref{dessin5} bellow and the family of 1-cycles $\gamma(t)$, represented on the $x$-plane by a big loop surrounding the four roots of $f(x)+t$,
on the family of elliptic curves $\Gamma_t$. Further,  
consider the complete elliptic integral $I(t)= \int_{\gamma(t)} xy dx$. A simple computation shows that the associated 0-dimensional Abelian integral is $\int_{\delta(z)}x^2$, where $\delta(z)=x_1(z)+x_2(z)+x_3(z)+x_4(z)$, $f(x_i(z))\equiv z$. On the other hand
 $$
 I'(t)= \int_{\gamma(t)} \frac{xdx}{2y}
 $$
 is a complete elliptic integral of third kind, and $\gamma(t)$ is homologous to a small loop around one of the two ''infinite" point of the affine curve $\Gamma_t$. The conclusion is that $I'(t)$ is a residue, in fact a non zero constant. The Abelian integral $I(t)$ is therefore linear in $t$. This example shows that the claim of Theorem \ref{main3} can not be improved.

\vskip 0.2cm
Note that the zero-cycles $\delta_i(z)$ 
are by no means unique, they depend on the mutual position of the segments $[c_0,c_i]$. If all the zero-cycles $\delta_i(z)$ are in the orbit of a given cycle $\delta_{i_0}$, obtained after a continuation with respect to $z$, then the vanishing of $\varphi_{i_0}$ implies the vanishing of all the $\varphi_{i}$, and hence of all the moments. Finally, the orbit of a given $\delta_{i_0}$ may contain other cycles, more suitable for our purposes. In the next subsection we propose an alternative construction of such a cycle, by using a residue calculus.
As we shall see, this will be more natural for the applications.

\subsection{The Cauchy integral related to $I$}
In this section we give an alternative computation of a convenient  necessary condition for the  identical vanishing of the Abelian integral $I(t)$, defined in (\ref{abelian}),  (\ref{abel}). Our result will hold under the additional assumption that there is a path along which the cycle $\gamma(t)$ vanishes. More precisely, 
let $\gamma(t)\subset \Gamma_t$ be a continuous family of closed continuous curves 
defined in a neighborhood of some regular value $t_0$ of $f$. Consider a path
\begin{equation}
\label{path}
[0,1] \rightarrow \C : s \mapsto t(s)
\end{equation}
such that $s(0)=t_0$, $s(1)=t_1$, $t(s)$ is a regular value of $f$  for $0\leq s < 1$, and $t_1$ is a singular value of $f$.
We shall say that the continuous family of closed loops $\gamma(t)$ vanishes along the path (\ref{path}) if it can be extended to a continuous family of loops along this path such that $\gamma(t_1)$ is homologous to zero on the singular affine curve $\Gamma_{t_1}$. This implies in particular that $I(t_1)=0$ as well as that the corresponding zero-dimensional Abelian integral vanishes at $t_1$.

Without loss of generality we suppose that 
$y$ restricted to $\gamma(t_0)$ does not vanish. 
Then, for all $(t,z)$ such that  $|z|$ and  $|t-t_0|$ are sufficiently small, the Cauchy type integral
\begin{equation}
\label{cauchy}
J_t(z)= \int_{\gamma(t)} \frac{k(x)y }{y^2-z} dx, \quad z\sim 0
\end{equation}
is well defined and analytic in $t,z$. 
 \begin{figure}[htbp]
\begin{center}
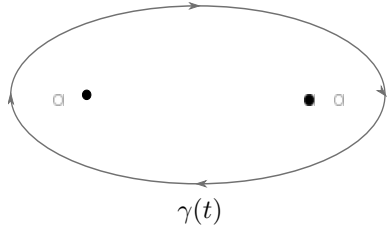
\caption{The definition of the Cauchy type integral $J_t(z)$.}
\label{dessin1}
\end{center}
\end{figure}
The definition of $J_t(z)$ is illustrated on fig.\ref{dessin1}, where a closed loop $\gamma(t)$ projected on the $x$-plane is shown, which makes one turn around two roots of the polynomial $f(x)+t$. The roots of $f(x)+t$ are represented by small black bullit circles, while the roots of $f(x)+t-z$ by small empty circles. Note that in section \ref{section2} we supposed that $z\sim \infty$, while
in this section that $z \sim 0$, and this is essential for what follows.

Since 
$
I'(t) = \frac{1}{2}  J_t(0) 
$, $J_t(z)$ is a deformation of the Abelian integral $I'(t)$. At the same time, for a fixed $t$,  $J_t(z)$ is a generating function of the moments
$I^{k}(t)$, $k\geq 1$, in the sense of section \ref{section2} and
$$
J_t(z)=J(z),
$$
where the Riemann surface $R=\Gamma_t$ depends on the parameter $t$. For $(x,y)\in \gamma(t)$ and $|z|$ sufficiently small the series
$$
\sum_{k=0}^\infty (\frac{z}{y^2})^k
$$
converges uniformly and hence
\begin{eqnarray*}
J_t(z)&=& \int_{\gamma(t)} \frac{k(x)y }{y^2(1-\frac{z}{y^2})} dx \\
&=& \int_{\gamma(t)} \frac{k(x)}{y }dx + z\int_{\gamma_0} \frac{k(x)}{y^3 }dx+ z^2\int_{\gamma_0} \frac{k(x)}{y^5 }dx + \dots
\end{eqnarray*}
Taking into consideration (\ref{ik1})
we conclude
\begin{proposition}
\label{41}
For every regular value $t$ of $f$ the equalities
$$ (^{-1/2}_{\mbox{  $k$}}) \frac{d^k}{dz^k}J_t(0)= 2 I^{(k+1)}(t), k= 0,1,2 \dots$$ hold.
\end{proposition}

The above Proposition implies the following corollary.
\begin{corollary}
The Abelian integral $I'(t)$ vanishes identically,  if and only if the Cauchy type integral
$J_{t}(z) $ vanishes identically. 
\end{corollary}

The main advantage of using $J_t(z)$ instead of $I'(t)$ is the possibility to extend it analytically with respect to  $z$.
The result is a function algebraic in $z$.
\begin{proposition}[\cite{pry04}]
For every fixed regular value $t$ the Cauchy type integral $J_t(z)$ extends to an algebraic function in $z$ with singularities at $z=0$ and at the critical values of $f$.
\end{proposition}
Indeed, for a fixed regular $t$, $J_t(z)$ allows for an analytic continuation along any path which does not contain critical values of $f-t$ or the value $z=0$. In a neighborhood of a critical value of $f-t$ or at $z=0$, the Cauchy theorem implies that, up to an addition of a holomorphic function, $J_t(z)$ is a linear combination of residues of $ \frac{k(x)y }{y^2-z} dx$ at the roots $f^{-1}_i(z-t)$ of $f(x)+t-z$. 
Thus, $J_t(z)$ is a function of moderate growth in $z$ with a finite number of branches, and hence is algebraic in $z$. \qed
\vskip 0.2cm

Our next goal is to extend analytically $J_t(z)$ in a neighborhood of $(t_1,0)$ under the condition that $t_1$ is a critical value of $f$. To simplify the notation put $t_1=0$, $f(0)=0$.  Consider the domain 
$$D_\delta= \{(t,z): |t|<\delta, |z|<\delta, t\neq z, t\neq 0, z\neq 0 \}$$
 and assume that $\delta>0$ is so small that $t=0$ is the only critical value of $f$ in the disc $\{t:  |t|<2\delta \}$. 
 Take some $(t,0)\in D_\delta$ and consider the germ of the analytic function $J=J_t(z)$ in a neighborhood of this point.
  \begin{figure}[htbp]
\begin{center}
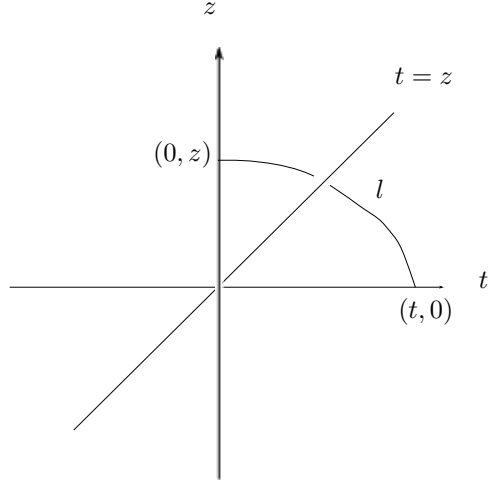
\caption{The domain $D_\delta$.}
\label{dessin2}
\end{center}
\end{figure}

\begin{figure}[htbp]
\begin{center}
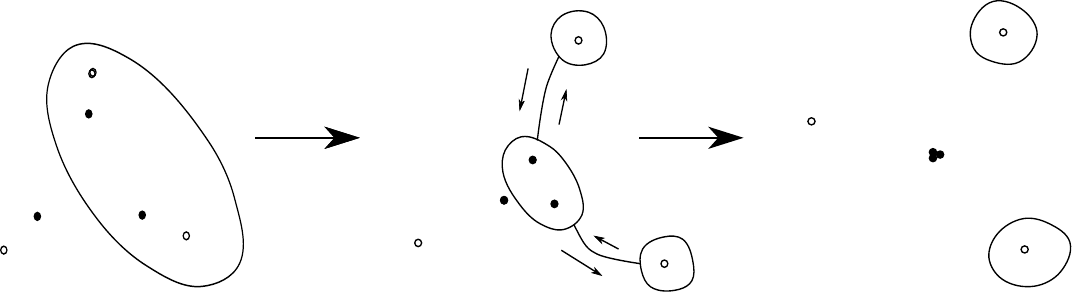
\caption{Computing the limit of $\gamma(t)$ at a singular value.}
\label{dessin4}
\end{center}
\end{figure}

\begin{figure}[htbp]
\begin{center}
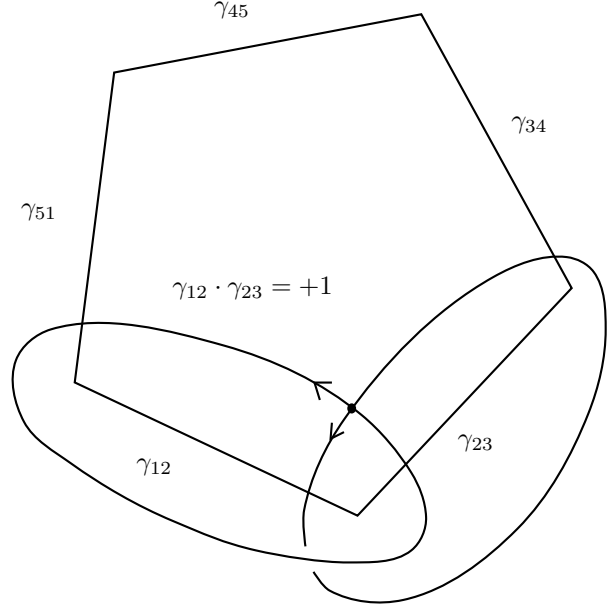
\caption{The standard basis of vanishing cycles of the the singularity $y^2+x^5$.}
\label{dessin3}
\end{center}
\end{figure}

\begin{proposition}
The germ of $J=J_t(z)$ at $(t,0)\in D_\delta$ allows for an analytic continuation along any path starting at $(t,0)$ and contained in 
$D_\delta$.
\end{proposition} 
Indeed,  the affine curve $\Gamma_t$ is regular, provided that $t\neq 0$, and the differential $ \frac{k(x)y }{y^2-z} dx$ has simple poles if and only if $z \neq t, z\neq 0$. Therefore the closed curve $\gamma(t)$ can be deformed in  a way to avoid these simple poles. \qed
\vskip 0.2cm
Although the function $J_t(z)$ might be not analytic in $t$ near the line $t=0$, it has a finite limit there which we compute next. For this purpose,
let $l  $ be a closed smooth path connecting the point $(t,0)$ to $(0,z)$, $t, z \neq 0$, and contained in $D_\delta$ (except the ends), see fig. \ref{dessin2}. Suppose that  the homology class of the limiting loop $\gamma(0)\subset \Gamma_0$ is zero and hence is a linear combination of vanishing cycles.  
\begin{theorem} 
\label{main4}
If $\gamma(0)\subset \Gamma_0$ is homologous to zero,
 then the limiting value of $J_t(z)$ at $(0,z)$ along $l$ is a zero-dimensional Abelian integral
$$
J_0(z) = 2\pi \sqrt{-z} \frac{d}{dz}\int_{\delta(z)} K(x)
$$
where $K(x)$ is a primitive of $k(x)$, $\delta(z) =  \sum_i n_i f^{-1}_i(z)$, $f^{-1}_i(z)$ are the roots of the polynomial $f(x)-z$, and the numbers $n_i$ depend only on the homology class represented by the loop $\gamma(0)$ in $H_1(\check{\Gamma}_0,\mathbb{Z})$, 
$\check{\Gamma}_0= \{ (x,y): y^2=f(x), f(x)\neq z \}$.
\end{theorem}
\begin{corollary}
\label{cormain}
If $I(t)= \int_{\gamma(t)} \omega \equiv 0$ then $\int_{\delta(z)} K \equiv 0$
\end{corollary}
\begin{corollary}
\label{corcycl}
According to Proposition \ref{41}, if $I'(t)=0$ for some regular $t$, then the multiplicity of this zero is the same as the multiplicity of $J_t(z)$ with respect to $z$ at $z=0$. In the particular case where $t=0$ is a Morse critical point, the Abelian integral $I'$ is analytic at $t=0$, and the multiplicity of the zero of $I'$ at $t=0$ is just the multiplicity of the zero of the analytic function  $J_0(z)$ at $z=0$. Thus, the multiplicity of the one-dimensional Abelian integral at a Morse critical point equals essentially  the corresponding multiplicity of the 0-dimensional Abelian integral.
\end{corollary}
\noindent{\it Proof of Theorem \ref{main4}.}
We can deform the loop $\gamma(t)$ along the interior of the path $l$ in a way to avoid the poles of $\frac{k(x)y }{y^2-z} dx$.
 Taking the limit $t\rightarrow 0$ along $l$ we obtain that $\gamma(0)$ is homologous to a sum of closed loops around the poles of $\frac{k(x)y }{y^2-z} dx$, as it is shown on fig. \ref{dessin4}.
Therefore,
\begin{eqnarray*}
J_0(z) & = & 2\pi \sqrt{-1} \sum_i n_i Res_{f^{-1}_i(z)} \frac{k(x) y}{f(x)-z} dx \\
 &= & 2\pi \sqrt{-z} \sum_i \pm n_i  \frac{k(f^{-1}_i(z))}{f'(f^{-1}_i(z))} \\
& = & 2\pi \sqrt{-z} \int_{\delta(z)} \frac{k(x)}{f'(x)}\\
& = & 2\pi \sqrt{-z} \frac{d}{dz}\int_{\delta(z)} K(x) .\ \ \ \Box
\end{eqnarray*}
{\bf Computation of the reduced 0-cycle $\delta(z)$.}
For simplicity, suppose that  $\gamma(t)$ vanishes as $t$ tends to $0$ at the origin $(0,0)$.
Thus $\gamma(t)$ is a linear combination of cycles vanishing at $(0,0)$. The  standard basis of such cycles can be described as follows. Let $f^{-1}_i(t)$, $i=1,2,\dots,n$ be the roots of the polynomial $f(x)+t$
which tend to $0$ as $t$ tends to $0$, ordered cyclically with respect to the monodromy action.
We denote by 
$\gamma_{ij}(t)\subset \Gamma_t$  a simple closed loop which is projected to the segment
$[f^{-1}_i(t),f^{-1}_j(t)]$. The loops $\gamma_{i,i+1}(t), i=1,2, \dots,n-1$ form a basis of the local homology group of the Milnor fiber of $y^2-f(x)$. We fix the orientations of these cycles by the convention
$$
\gamma_{i,i+1}\cdot \gamma_{i+1,i+2}(t)=1 .
$$
It is easy to check that then
$$
\gamma_{i,i+1}(t)+\gamma_{i+1,i+2}(t)= \gamma_{i,i+2}(t),
$$
where the orientation of $\gamma_{i,i+2}(t)$ is appropriately chosen.
Therefore the orientations of the remaining cycles can be chosen to satisfy
\begin{equation}
\label{ij}
\gamma_{ij} \circ \gamma_{jk} = +1, \gamma_{ij} + \gamma_{jk}= \gamma_{ik}, \ \ \forall i<j<k .
\end{equation}
As a by product we have also
$$
\gamma_{1,2}+ \gamma_{2,3}+ \dots + \gamma_{n,1} = 0 .
$$
Obviously this fixes the orientation of all cycles $\gamma_{ij}$ up to simultaneous multiplication by $-1$, which have no incidence on the result claimed in Corollary \ref{cormain}. The standard basis of vanishing cycles of the the singularity $y^2+x^5$ is shown on fig.\ref{dessin3}.
We shall construct an isomorphism
\begin{eqnarray*}
H_1(\Gamma_t,\Z)&\rightarrow&\tilde{H}_0(f^{-1}(z), \Z)\\
\gamma(t) &\mapsto &\delta(z)
\end{eqnarray*}
 having the property announced in Corollary \ref{cormain}. According to the proof of Theorem \ref{main4} this should be a linear map which associates to the one-cycle $\gamma_{ij}(t)$  the reduced 0-cycle (see also section \ref{section2})
$\delta_{ij}(t)=  \pm (f^{-1}_i(t)- f^{-1}_j(t))$ and should be therefore compatible to the relations
$$
 \gamma_{ij} + \gamma_{jk}= \gamma_{ik},  \quad \delta_{ij} + \delta_{jk}= \delta_{ik}.
$$
It follows that the orientation of the 0-cycles $\delta_{ij}(t)$ can be fixed as
$$\delta_{ij}(t)=  f^{-1}_i(t)- f^{-1}_j(t),\ \ \ \forall i<j .$$
Note that the above isomorphism is not compatible to the monodromy action.

In conclusion, if
$$
\gamma(t) = \sum  n_{ij} \gamma_{ij}(t)
$$
and 
 $I(t)= \int_{\gamma(t)} \omega \equiv 0$ then $\int_{\delta(z)} K \equiv 0$
where
$$\delta(z)= \sum  n_{ij} (f^{-1}_i(z) - f^{-1}_j(z)).$$
\subsection{Hyperelliptic Abelian integrals along ovals.} 
Let $f(x)\in \R[x]$ be an arbitrary  non-linear real polynomial.  Consider a  family of ovals 
$\{\gamma(t)\}_t$ 
$$
\gamma(t) \subset \{(x,y)\in \R^2: y^2-f(x)=t \} , t\in \R
$$
depending continuously on the real parameter $t$.
Each oval $\gamma(t)$ can be parameterized as
$$
y=\pm\sqrt{f(x)+t}, \quad x_1(t)\leq x \leq x_2(t)
$$
where $x_1(t)< x_2(t)$ are  two real roots of $f(x)+t$. The purpose of this last section is to solve, by making use of Theorem \ref{main3} and Theorem \ref{main4}, the following problem: under what conditions the Abelian integral (\ref{abel})
$$
I(t)=\int_{\gamma(t)} k(x)y dx = 2 \int_{x_1(t)}^{x_2(t)} k(x)y dx
$$
is identically zero? 

\begin{theorem}
\label{exth}
The integral $I(t)$ vanishes identically if and only if there exists a  polynomial 
$ r\in \R[x]$, such that both $f$ and $K=\int k$ are polynomials in $r$, and $ r(x_1(t))\equiv r(x_2(t))$.\\
%$f$ and $K=\int k$ have a right compositional polynomial factor, identifying $x_1(t)$ and $x_2(t)$.
\end{theorem}
\pr
First of all, note that if $K$ and $f$ have a right compositional factor identifying $x_1(t)$ and $ x_2(t)$, then the Abelian integral
$\int_{\gamma(t)} k(x)y dx$ is a pull back of an integral along a cycle homologous to zero, and hence vanishes identically. 

Suppose further that $I(t)$ vanishes identically. It is enough to show that this implies
$K(x_1(t)) \equiv K(x_2(t))$ since in this case by Proposition \ref{ecc} (or by the L\"{u}roth theorem) $f$ and $K$ will have a right compositional factor identifying $x_1(t)$ and $x_2(t)$. If there exists a path on the complex $t$-plane along which the cycle $\gamma(t)$ vanishes, then Theorem \ref{main4} applies and we conclude that $K(x_1(t)) \equiv K(x_2(t))$.

As an example, consider
a real polynomial $f$  of degree $n\geq 2k$, $f = -x^{2k} + \dots$. Let 
 $x_1(t)< x_2(t)$ be the two real roots of $f(x)+t$ which tend to $0$ as $t$ tends to zero and $\{ \gamma(t) \}$ be the  continuous family of ovals 
 vanishing at the origin as $t$ tends to zero. 
 $$\gamma(t)\subset \{(x,y)\in \mathbb{R}^2: y^2+x^{2k}+ \dots = t\}$$
Then Theorem \ref{main4} applies  and hence the result of Theorem \ref{exth} follows. In the Morse case ($k=1$), this has been proved by  Christopher and Mardesic \cite{krma10}.

The condition that $\gamma(t)$ vanishes along a suitable path is essential, and holds for arbitrary real polynomials of degree  four or five, see for instance \cite[section 3.1]{gail05}, where the case $f(x)=(x^2-1)^2$ is studied. We do not know whether this condition is fulfilled for arbitrary polynomial $f$ and family of ovals $\gamma(t)$. See Fig. 8 (continuous families of ovals). However, using Theorem \ref{main3} instead of Theorem \ref{main4}
we can
prove the theorem in its full generality.
   
\begin{figure}[htbp]
\begin{center}
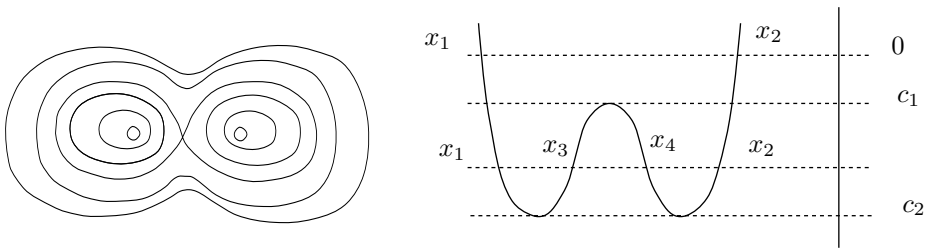
\caption{The continuous families of ovals of $y^2+(x^2/2-1)^2$ and the graph of $(x^2/2-1)^2-1$.}
\label{dessin5}
\end{center}
\end{figure}

Indeed, let $f$ be an arbitrary real polynomial of degree $n>1$ and $I(t)$ be an identically vanishing Abelian integral as before. Let us apply Theorem \ref{main3}. For this purpose, let us fix a regular real value $t$ of $f$, and consider the moment problem associated to the oval $\gamma(t)$ on the Riemann surface $\Gamma_t$. Following the method described in section \ref{section2} we have to consider a constellation $\lambda_f \subset \P^1$ and to deform the image of $\gamma(t)$ under $f+t$ on $\lambda_f$. The closed loop $\gamma(t)$ being an oval, its image is just a real interval connecting $0$ to a critical value  $c_k$ of $f+t$. Suppose for instance that $0>c_1>c_2>\dots > c_k$ are the remaining critical values of $f+t$ contained in $[c_k,0]$. We have therefore
$$
[c_k,0]= [c_k,c_{k-1}]\cup\dots\cup [c_1,0]
$$
which, without loss of generality, will be used on the place of the constellation $\lambda_f$. 
To each segment $[c_{i-1},c_i]$ we associate a 0-cycle $\delta_i$ and $I(t)$ is a rational function if and only if $\int_{\delta_i(z)} K \equiv 0$. 
\begin{quote}
Example. \emph{The critical values of the polynomial $(x^2/2-1)^2-1$ are $-1$ and $-3/4$. The relevant constellation associated to the exterior family of ovals shown on  fig.\ref{dessin5} is $[-1,-3/4]\cup[-3/4,0]$. To the segment $[-3/4,0]$ we associate the 0-cycle $\delta_1(z)=x_1(z)-x_2(z)$ and to the segment $[-1,-3/4]$ the 0-cycle
$\delta_2(z)=x_1(z)-x_3(z)+x_4(z)-x_2(z)$.}
\end{quote}
Fortunately in general we do not need to compute all of $\delta_i$. We note that the image of $[x_1(t),x_2(t)]$ is a closed curve covering $[c_k,0]$. The pre-image of each point  $z\in(0,c_1)$ consists of two points $x_1(z)$ and $x_2(z)$ (roots of $f(x)+t-z$) and hence $\delta_1(z)= x_1(z) - x_2(z)$. We conclude that $K(x_1(t)) \equiv K(x_2(t))$ which completes the proof of Theorem\ref{exth}. \qed

\bibliographystyle{amsplain}
%\bibliography{../BIBfiles/bibliography}

%\end{document}

\end{document}